\newcounter{version}
\documentclass[12pt]{amsart}

\ifnum\theversion>0
  
\fi
\numberwithin{equation}{section}

\newtheorem{thm}{Theorem}[section]
\newtheorem{cor}[thm]{Corollary}

\newtheorem{lem}[thm]{Lemma}

\theoremstyle{definition}

\theoremstyle{remark}
\newtheorem*{rem}{Remark} 

\newcommand{\C}{{\mathbb C}}
\newcommand{\R}{{\mathbb R}}
\newcommand{\Proj}{{\bf P}}
\newcommand{\Hypb}{{\bf H}}
\newcommand{\dbar}{\overline{\partial}}
\newcommand{\Fm}{\Lambda}
\newcommand{\abs}[1]{\lvert#1\rvert}

\newcommand{\Hom}{\qopname\relax o{Hom}}
\newcommand{\vecsp}[1]{\langle #1\rangle}
\newcommand{\Un}{\qopname\relax o{U}}
\newcommand{\PU}{\qopname\relax o{PU}}
\newcommand{\im}{\sqrt{-1}}
\newcommand{\Hl}{\Hypb^l}
\newcommand{\cphi}{\Phi}
\newcommand{\cpsi}{\Psi}
\newcommand{\tildeU}{U'}

\title{Twistor spaces of non-flat Bochner-K\"ahler manifolds}
\author{Yoshinari Inoue}
\address{Department of Mathematics\\
        Faculty of Science\\ Kyoto University\\
        Kyoto 606-01, Japan}
\email{inoue@kusm.kyoto-u.ac.jp}
\usepackage{amsmath,amsthm}

\begin{document}
\maketitle

\section{Introduction}

It was shown by Kamishima in \cite{[Ka]} that
a Bochner-K\"ahler manifold of dimension $n>1$
 is locally isomorphic to one of
the following model spaces:
\begin{enumerate}
\item
a flat geometry $(\C^n\times (\Un(n)\times \R_+), \C^n)$,
\item
a non-flat geometry
 $(\PU(1,l)\times \PU(m+1), \Hl\times\Proj^m)$.
\end{enumerate}
This can be compared to the fact that a conformally flat
manifold of dimension $n>2$ is locally conformal to a region
of the sphere of the same dimension.
In twistor theory, it is well-known that an even dimensional
conformally flat manifold has an
integrable twistor space
(\cite{[O.R]}, \cite{[M]}, \cite{[I1]}, \cite{[I3]}).
It is interesting, as an analogy, that
a Bochner-K\"ahler manifold has integrable twistor spaces
defined by O'Brian and Rawnsley in \cite{[O.R]}.
Since it was proved by calculating
the Nijenhuis tensors, the proof does not  give holomorphic
coordinates. Hence it would be natural to ask to
``give a system of local coordinates of the
twistor spaces of the above model spaces''.
Furthermore, to study a generalization of the Penrose transform,
it is necessary to ``construct the moduli space of relative deformations
of fibers'', which we call the {\it complexification}.

In the present paper, we give answers to the above two questions.

Since the first space is flat, its twistor space as a Riemannian manifold
has already an integrable complex structure.
Hence the problems can be solved easily by restricting the Penrose
diagram as a Riemannian manifold.

In the second case,
we give local coordinates by constructing a local embedding
to a product of complex projective spaces
(Theorem \ref{hol} and Theorem \ref{embd}).
The complexification of $\Hypb^l\times\Proj^m$
is $\Hl\times {\bar\Hypb}^l \times
	\left(\Proj^m\times \bar{\Proj}^m \setminus F(1,m;m+1)\right)$
(Corollary \ref{cmlx}).

Let us explain briefly the contents of this paper.
In \S2, we review the definition of the twistor spaces of Hermitian
manifolds. The answers in the flat case are explained shortly.
In \S3, we construct a system of local coordinates
of the twistor spaces in the second case.
In \S4, we give an explicit description of the complexification
of the model space, by extending maps in the previous section.
In \S5, we give some examples.

\section{The twistor spaces of Hermitian manifolds}

In this section, we review the definition of the twistor spaces
of a general Hermitian manifold $M$.
Let $n$ be the complex dimension of $M$ and
 $k$  an integer between $0$ and $n$.
  The $k$-th twistor space $Z_k(M)$ of $M$
is the total space of the Grassmannian bundle $G_k(T^{(1,0)}M)$
equipped with the almost complex structure defined as follows.
By the Hermitian metric on $M$, we can define a connection on its
 frame bundle.
Since $G_k(M)$ is an associated bundle of the frame bundle,
the connection decomposes
the real tangent bundle of $Z_k(M)$ as
\begin{equation*}
	T Z_k(M) =  T_V Z_k(M) \oplus T_H Z_k(M).
\end{equation*}
Since its fibers have natural holomorphic structures,
we have a complex structure on the first component.
For $x\in M$ and $L\in G_k(T^{(1,0)}_xM)$,
decompose the second component as
$$
	T_H Z_k(M)_{(x,L)} = T^{(1,0)}_x M = L \oplus L^{\perp}.
 $$
This decomposition is compatible with the original complex structure,
since $L$ is a complex subspace of $T^{(1,0)}_xM$.
We define a new complex structure on ${T_H} Z_k(M)_{(x,L)}$
by replacing the complex structure
on $L$ with its conjugate structure.

We can write
explicitly the space of horizontal $(1,0)$-forms,
by using the Pl\"ucker coordinates
with notation of multi-indices similar to
that in \cite{[I2]} and \cite{[I3]}.
Let $(e_1,\dots,e_n)$ be a local orthonormal basis of
the space of $(1,0)$-vectors of $M$.
For $I_0=(i_1\dots i_k)$ with $i_1 < \dots < i_k$,
let $z^{I_0}$ be the coefficient of
$e_{I_0} = e_{i_1}\wedge\dots\wedge e_{i_k}$.
A multi-index $I$ is a sequence of elements of $\{1,\dots,n\}$.
The composition of the sequences $I$ and $J$ is written as $IJ$.
We put relations between $z^I$'s as:
\begin{align*}
	z^{I_1iiI_2} &= - z^{I_1I_2},\\
	z^{I_1ijI_2} &= - z^{I_1jiI_2}, \quad i\not=j.
\end{align*}
Then, for any multi-index $I$, there is a unique subsequence
$I_0$ of $(1,\dots,n)$ such that $z^I = \pm z^{I_0}$.
We denote $\abs{I}$ to be the length of $I_0$.
Then the space of horizontal $(1,0)$-forms
of the twistor space is spanned by the following forms:
\begin{equation}\label{ACS}
  \begin{split}
	\sum_{j\not\in I} z^{jI} \bar{e}^j,&\quad	\abs{I}=k-1,\\
	\sum_{i\in I}  z^{iI} e^i,&\quad \abs{I}=k+1,
  \end{split}
\end{equation}
where $(e^1,\dots,e^n)$ is the dual basis.

 The integrability condition of the almost complex structure
is given by O'Brian and Rawnsley
in \cite{[O.R]}. A Bochner-K\"ahler manifold,
which we will study in the following sections, is an example
of such manifolds.

Henceforth, we assume that the almost complex structure of $Z_k(M)$
is integrable.

We can define the complexification of
an even dimensional conformally flat manifold
by considering the relative deformations of fibers of its twistor space.
Since Grassmannian manifolds are rigid, we can
also define the {\it complexification} of $M$ in a similar way.

For $(A,B)\in \Un(k)\times\Un(n-k)$, define
\begin{align*}
	\rho(A,B): \C^k\oplus \C^{n-k} &\rightarrow \C^k\oplus\C^{n-k}\\
		    (v_1,v_2) &\mapsto (A^*v_1, Bv_2).
\end{align*}
Then, for $x\in M$, the normal bundle $N_x$
of the fiber $Z_k(M)_x$ is written as
$$
	N_{x} \simeq N = \Un(n) \times_{\rho} (\C^k\oplus\C^{n-k}).
 $$
By the theorem of Bott-Borel-Weil-Kostant, we have:
$$
	H^i(Z_k(M)_x, {\mathcal O}(N_x)) =
	  \begin{cases}
	    \Fm^{(1,0)}_x M \oplus T^{(1,0)}_x M,&\quad \text{if $i=0$,}\\
	    0,& \quad \text{if $i>0$.}
	  \end{cases}
 $$
Hence, by the relative deformation theory of Kodaira (\cite{[Ko]}),
the moduli space
$$
	M_\C = \{ Z\subset Z_k(M)\mid \text{$Z$ is a deformation of
	 a fiber.} \}
 $$
is a $2n$-dimensional complex manifold, which includes $M$
as the space of {\it real} fibers.

 As an example, we describe the situation shortly when $M=\C^n$ with
flat metric.

 Let $X= \C^n\times \bar{\C}^n	= H^0(G_{k,n}, \mathcal{O}(N))$.
Let $Y$ be the submanifold of
$N\times X$ defined by
$$
	Y = \{ (f(z),f) \mid f\in
	  X, z\in G_{k,n} \}.
 $$
We have a canonical double fibration:
\begin{equation*}
\setlength{\unitlength}{0.5mm}%
  \begin{picture}(50,30)
    \put(20,20){$Y$}
    \put(0,0){$N$}
    \put(40,0){$X.$}
    \put(18,18){\vector(-1,-1){10}}
    \put(27,18){\vector(1,-1){10}}
    \put(7,16){$\scriptstyle p_1$}
    \put(34,16){$\scriptstyle p_2$}
  \end{picture}
\end{equation*}
Then we can show easily that, for each $w\in N$, there is
a unique element $v\in \C^n$ such that
$w\in p_1\circ p_2^{-1}(v,\bar{v})$.
By the projection $w\mapsto v$, we can identify $N$ with
$Z_k(\C^n)$.
The {\it complexification} of $\C^n$ is $X$, namely,
for $x\in X$,  the corresponding
submanifold is written as $p_1\circ p_2^{-1}(x)$.
\begin{rem}
Since the twistor space $Z(\C^n)$ as a Riemannian manifold
satisfies the integrability condition,
this diagram can be obtained by restricting
the Penrose double fibration as a Riemannian manifold,
which is described in detail by Murray in \cite{[M]}.
\end{rem}

\section{The model spaces of non-flat Bochner-K\"ahler manifolds
and their twistor spaces}

Let $V$ be a complex vector space of dimension $l+1$  having
the indefinite metric:
$$
	\begin{pmatrix}
	  1 & 0\\
	  0 & -I_l
	\end{pmatrix}.
 $$
Then we can consider the $l$-dimensional complex hyperbolic space
$\Hl = \Hl(V)$ as an open submanifold of $\Proj^l = \Proj(V)$ by
$$
	\Hl = \{x = (1: x_1 : \dots : x_l)\in \Proj^l \mid
	  \abs{x_1}^2 + \dots + \abs{x_l}^2 < 1 \}.
 $$
Let $W$ be a Hermitian vector space of dimension $m+1$
and  $\Proj^m=\Proj(W)$ its projective space.
Let $y = (1:y_1:\dots:y_m)$ denote its local affine coordinates
with respect to an orthonormal basis of $W$.

We consider $\Hl$ and $\Proj^m$	 as K\"ahler manifolds
by the well-known K\"ahler forms
\begin{align*}
	& -\im \partial \bar{\partial}
	  \log \left(1- \sideset{}{_i}\sum \abs{x^i}^2\right),\\
	& \im \partial \bar{\partial}
	  \log\left(1+\sideset{}{_j}\sum \abs{y^j}^2\right),
\end{align*}
respectively.
Then $M=\Hl\times \Proj^m$ with the product metric is
an example of Bochner-K\"ahler manifolds.
Hence
the almost complex structure of $Z_k(M)$ is integrable.
By the definition of the metric,
the group $\Un(1,l) \times \Un(m+1)$ acts on $M$ isometrically.

Let $x_0=(1,0,\dots,0)$ and $y_0=(1,0,\dots,0)$.
By identifying a tangent vector with an element of
$x_0^\perp\oplus y_0^\perp$,
we can consider it as an element of $V\oplus W$.
Let $z^{IJ'}$ be the Pl\"ucker coordinate,
where $I$ (resp. $J$) is a multi-index
with respect to $x_1,\dots,x_l$ (resp. $y_1,\dots,y_m$).
Note that the local frame of the tangent bundle that
we choose is orthonormal only at $(x_0,y_0)$.
Hence \eqref{ACS} is valid only at points on the fiber over $(x_0,y_0)$.

Let $L$ be a complex vector space with metric.
For a vector $x\in L$, we regard $\bar{x}$ as
 an element of $L^*$ by the canonical
isomorphism $\bar{L} \simeq L^*$.
Then, by using the decomposition:
$$
	\Fm^{k+1}(V\oplus W) = \bigoplus_{a+b=k+1}\Fm^aV\otimes\Fm^bW,
 $$
we define
\begin{align*}
  \phi_{a,b} &= \left.
  x\wedge\left(1-\frac{y\wedge i(\overline{y})}{(y,y)}\right) z +
    y\wedge\left(1-\frac{x\wedge i(\overline{x})}{(x,x)}\right) z
  \right|_{\Fm^{a} V\otimes \Fm^{b} W}\\
  &= \sum_{\abs{I}=a, \abs{J}=b} \phi_{I,J} e_I\wedge e_{J'},
\end{align*}
where $i$ denotes the interior multiplication.
In the following, we assume that multi-indices $I,J,\dots$
do not include the index $0$.
Then we have
\begin{equation}\label{a}
\begin{split}
	\phi_{0I,0J} &= (-1)^a
	  \left(
	    \sum_{i\not\in I}\frac{\overline{x}^i z^{iIJ'}}{(x,x)} +
	    \sum_{j\not\in J}\frac{\overline{y}^j z^{j'IJ'}}{(y,y)}
	  \right),\\
	\phi_{I,0J} &= (-1)^a z^{IJ'} - \sum_{i\in I} x^i \phi_{0iI,0J},\\
	\phi_{0I,J} &= z^{IJ'} - \sum_{j\in J} y^j \phi_{0I,0jJ},\\
	\phi_{I,J} &= - \sum_{j\in J} y^j z^{j'IJ'}
	  - \sum_{i\in I}x^i \phi_{0iI,J},
\end{split}
\end{equation}
where $\bar{x}^i$ is the complex conjugate of $x^i$.

Let $p_1$ and $p_2$ be the projections:
\begin{align*}
	p_1: \Hl\times \Proj^m &\rightarrow \Hl,\\
	p_2: \Hl\times \Proj^m &\rightarrow \Proj^m.
\end{align*}
Let $L_{\Hl}$ and $L_{\Proj^m}$ be the tautological line bundles.
Let $L_{Z_k(M)}$ denote
the tautological line bundle
as a Grassmannian bundle.
Then $\phi_{a,b}$ is canonically extended to a map
$$
	\phi_{a,b}: L_{Z_k(M)}\otimes p_1^*(L_{\Hl}^a)
			\otimes p_2^*(L_{\Proj^m}^b)
		\rightarrow
	\Fm^a V\otimes \Fm^b W,
 $$
which is equivariant under the
action of $\Un(1,l)\times \Un(m+1)$.
The space on the left hand side has an almost complex structure
induced by that of $Z_k(M)$ and the canonical connection
as a line bundle. Then we have the following theorem.

\begin{thm}\label{hol}
The map $\phi_{a,b}$ is holomorphic.
\end{thm}

\begin{proof}
It suffices
to calculate that $\phi_{a,b}$ is $\dbar$-closed at points on the
fiber over $(x_0,y_0)$.
Then it is immediate by the definition
of the almost complex structure.
\end{proof}

\begin{rem}
This theorem shows that the line bundle
$L_{Z_k(M)}\otimes p_1^*(L_{\Hl}^a)
		\otimes p_2^*(L_{\Proj^m}^b)$
has a natural holomorphic structure if $a+b=k+1$.
In fact, if both $l$ and $m$ are positive,
this is necessary,
which can be shown directly by computing the integrability condition.
\end{rem}

For non-negative integers $a$ and $b$ such that $a+b=k$, put
$$
	U_{a,b} = \{ \kappa \in Z_k(M) \mid
	  \phi_{a+1,b}(\kappa) \not= 0, \phi_{a,b+1}(\kappa) \not= 0 \}.
 $$
These open sets cover $Z_k(M)$. In fact, we have the following
lemma.
\begin{lem}
Let $(x,y,z)$ be a point of $Z_k(M)$ and
 $L$ the subspace of
$x^\perp\oplus y^\perp$ corresponding to $z$. Then
$(x,y,z)\in U_{a,k-a}$ if and only if
$\dim (L\cap x^\perp) \le a \le k - \dim(L\cap y^\perp)$.
\end{lem}

\begin{proof}
By \eqref{a}, we have
\begin{align}
	\phi_{0I,0J}(x_0, y_0, z) &= 0,\label{eq1}\\
	\phi_{I,0J}(x_0, y_0, z) &= (-1)^{\abs{I}} z^{IJ'},\label{eq2}\\
	\phi_{0I,J}(x_0, y_0, z) &= z^{IJ'},\label{eq3}\\
	\phi_{I,J}(x_0, y_0, z) &= 0\label{eq4}.
\end{align}
Hence, by the Pl\"ucker relations, we have
$$
	U_{a,b} = \{ (x,y,z)\in Z_k(M) \mid
	  L |_{\Fm^a x^{\perp}\otimes \Fm^b y^{\perp}} \not=0 \},
 $$
from which the lemma follows immediately.
\end{proof}

Now we define a map:
\begin{align*}
	\psi_{a,b}:
	 U_{a,b} &\rightarrow \Proj(\Fm^{a+1}V\otimes\Fm^{b}W)
		\times \Proj(\Fm^{a}V\otimes\Fm^{b+1}W)\\
	 \kappa &\mapsto ([\phi_{a+1,b}(\kappa)],[\phi_{a,b+1}(\kappa)]).
\end{align*}

\begin{thm}\label{embd}
The map $\psi_{a,b}$ is an embedding.
\end{thm}

\begin{proof}
Since the map is equivariant,
to prove the injectivity,
 it suffices to prove that
if $\psi_{a,b}(x,y,w) = \psi_{a,b}(x_0,y_0,z)$ then $x=x_0$,
$y=y_0$ and $w=z$.
Suppose $(y,y_0) = 0$. Then we can assume $y = (0,1,0,\dots,0)$.
We assume, for a while, that $J$ does not include the index $1$.
We compute
$$
	\phi_{0I,0J} = 0 = w^{I0'J'}, \quad \abs{I} = a, a-1.
 $$
Hence we have
\begin{gather*}
	\phi_{I,0J} = (-1)^{\abs{I}} z^{IJ'} =
	  - \sum_{i\in I} x^i \phi_{0iI,0J} = 0,\quad \abs{I} = a,\\
	\phi_{I,01J} = (-1)^{\abs{I}}z^{I1'J'} =
	- (-1)^{\abs{I}}w^{I0'J'} = 0,\quad \abs{I} = a.
\end{gather*}
These contradict with the assumption $z\in U_{a,k-a}$.

Hence we should have $(y,y_0) \not=0$.
Then, by \eqref{eq2} and \eqref{eq3},
there exists a constant $c$ such that
$w^{IJ'} = c z^{IJ'}$ for all $(I,J)$ with $\abs{I} = a-1, a, a+1$.
Hence we have $w=z$ by the Pl\"ucker relations.
By \eqref{eq1} and  \eqref{eq4}, the equations for $x$ and $y$ become
\begin{align*}
	\sum_{i\not\in I}
	    \frac{\overline{x}^i z^{iIJ'}}{(x,x)}
	+\sum_{j\not\in J}
	    \frac{\overline{y}^j z^{j'IJ'}}{(y,y)} &= 0,\\
	 \sum_{j\in J} y^j z^{j'IJ'} + \sum_{i\in I} x^i z^{iIJ'} &= 0.
\end{align*}
Let $x'=(x^1,\dots,x^l)$ and $y'=(y^1,\dots,y^m)$. Then these equations
 mean that there is a $k$-dimensional subspace
$L$ of $x_0^\perp\oplus y_0^\perp$ such that
$$
	(x',y') \in L \subset
	  (\frac{x'}{(x,x)}, \frac{y'}{(y,y)})^\perp.
 $$
Hence we have $x=x_0$ and $y=y_0$, concluding that the map
is injective.

 Injectivity of the derivative is also computed by
simple computations at points on the fiber over $(x_0, y_0)$.
\end{proof}

\section{The complexification of the model spaces}

In this section, we give an explicit description
of the {\it complexification} $M_\C$
of the model space $M = \Hl\times\Proj^m$.

Let $\xi, \zeta\in \Hl$ and
$\mu, \nu\in\Proj^m$ such that $(\mu,\nu)\not=0$.
We construct an embedding of
$G_{k,n} = G_k(V/\vecsp{\xi} \oplus W/\vecsp{\mu})$
in $Z_k(M)$ by
generalizing $\phi_{a,b}$	in the previous section.
Let $a$ and $b$ be non-negative integers such that $a+b=k+1$.
For $z\in G_k(V/\vecsp{\xi}\oplus W/\vecsp{\mu})$, we define:
$$
	\cphi_{a,b}(z) =
	  \left.
	  \xi\wedge
	    \left(
		1 - \frac{\mu\wedge i(\bar{\nu})}{(\mu,\nu)}
	    \right)z
	  + \mu\wedge
	    \left(
		1 - \frac{\xi\wedge i(\bar{\zeta})}{(\xi,\zeta)}
	    \right)z
	 \right|_{\Fm^aV\otimes\Fm^bW\,.}
 $$
For non-negative integers $a$ and $b$ satisfying $a+b=k$,  put
\begin{align*}
	\tildeU_{a,b} &= \{ z \in G_{k,n} \mid
	  \cphi_{a+1,b}(z) \not= 0, \cphi_{a,b+1}(z) \not= 0 \}\\
	&= \{ z \in G_{k,n} \mid
	  z |_{\Fm^a V/\vecsp{\xi}\otimes \Fm^b W/\vecsp{\mu}}
		 \not=0 \}.
\end{align*}
Then we define
\begin{align*}
	\cpsi_{a,b}(\xi,\zeta,\mu,\nu)=	\cpsi_{a,b}:
	 \tildeU_{a,b} &\rightarrow \Proj(\Fm^{a+1}V\otimes\Fm^{b}W)
		\times \Proj(\Fm^{a}V\otimes\Fm^{b+1}W)\\
	 z &\mapsto ([\cphi_{a+1,b}(z)],[\cphi_{a,b+1}(z)])
\end{align*}

\begin{thm}
Let $z\in G_{k,n}$.
There exists $\kappa\in Z_k(M)$ such that
$\psi_{a,b}(\kappa) = \cpsi_{a,b}(z)$
for all $(a,b)$ such that $z\in \tildeU_{a,b}$.
Hence $\cpsi_{a,b}$ induces an embedding of $G_{k,n}$ in $Z_k(M)$.
\end{thm}

\begin{proof}
Let $L$ be the subspace of $V/\vecsp{\xi}\oplus W/\vecsp{\mu}$
corresponding to $z$.
We can assume $\zeta = x_0$ and $\nu = y_0$.
By the isomorphisms
$x_0^{\perp} \simeq V/\vecsp{\xi}$ and
$y_0^{\perp} \simeq W/\vecsp{\mu}$
defined by the canonical projections,
we can consider $L$ as a subspace of $x_0^{\perp}\oplus y_0^{\perp}$.
We give a point $\kappa = (x, y, w)$ which satisfies
the property. We fix $a$ and $b$ for a while, and calculate
 the coefficients
of $\cphi_{a+1,b}(z)$ and $\cphi_{a,b+1}(z)$.
By a similar argument
in the proof of Theorem~\ref{embd}, we should have
$(y,y_0) \not= 0$. Thus, by \eqref{a},	we have
\begin{align*}
	\cphi_{0I,0J}(z) &=  0,\\
	\cphi_{I,0J}(z) &= (-1)^{\abs{I}}z^{IJ'},\\
	\cphi_{0I,J}(z) &= z^{IJ'},\\
	\cphi_{I,J}(z) &= - \sum_{j\in J} \mu^j z^{j'IJ'}
		   - \sum_{i\in I} \xi^i z^{iIJ'}.
\end{align*}
Hence, by the Pl\"ucker relations, we have $w=z$.
The coordinates of $x$ and $y$ satisfy the equation:
\begin{equation}\label{int}
	\sum_{i\not\in I} \frac{\bar{x}^i z^{iIJ'}}{(x,x)}
	  +
	\sum_{j\not\in J} \frac{\bar{y}^j z^{j'IJ'}}{(y,y)}
	  = 0.
\end{equation}
Define
\begin{align*}
	V_1 &= L \cap x_0^\perp,\\
	W_1 &= L \cap y_0^\perp.
\end{align*}
Let $V_2$ and $W_2$ be the subspaces such that
\begin{align*}
	x_0^\perp &= V_1 \oplus V_2\\
	y_0^\perp &= W_1 \oplus W_2
\end{align*}
which are orthogonal decompositions.
Then \eqref{int} means that
\begin{align*}
	x &\in \vecsp{x_0} + V_2,\\
	y &\in \vecsp{y_0} + W_2.
\end{align*}
Hence the problem is reduced to the case $V_1 = W_1 = 0$,
which implies that $\min(l,m)\ge k$.
By considering the canonical anti-holomorphic diffeomorphism:
\begin{align*}
	Z_k(M) &\rightarrow Z_{n-k}(M)\\
	L &\mapsto L^\perp,
\end{align*}
we can further assume that
\begin{gather*}
	\dim V = \dim W = k+1,\\
	z \in \bigcap_{0\le a \le k} \tildeU_{a,k-a}.
\end{gather*}
In this case, it is easily verified that the
equations for $x$ and $y$ are independent with respect to $a$.
Hence it suffices to find $x$ and $y$ when $a=0$.

An image of the map
$$
	\cpsi_{0,k}: \bigcap_{0\le a \le k} \tildeU_{a,k-a}
	  \rightarrow \Proj(\Hom(W, V))
 $$
is a projectivized invertible linear map.
Hence it suffices to find
$x\in V$ and $y\in W$ such that
\begin{align*}
	A(y) &= x,\\
	A(y^\perp) &= x^\perp,
\end{align*}
for an arbitrary invertible linear map $A: W\rightarrow V$.
Indeed, if we find them, we have
$$
	\psi_{0,k}\left (x,y, \vecsp{A(e_j) + e_j \mid j=1,\dots,k}\right)
	  = [A],
 $$
where $(e_1,\dots,e_k)$ is a basis of $y^\perp$.

We denote the indefinite metric on $V$ as $(,)_G$.
By the automorphism $A$, we define a Hermitian metric $(,)_H$ on $V$.
For a positive number $r$, let $S_r$ be the real hypersurface of $V$
defined as
$$
	S_r = \{ v\in V \mid (v,v)_H = r \}.
 $$
Let $H$ be the real hypersurface of $V$ defined as
$$
	H = \{ v\in V \mid (v,v)_G = 1 \}.
 $$
Let $r_0$ be the minimum number such that $S_r\cap H \not = \emptyset$.
Then we can show easily that
$S_{r_0}$ and $H$ intersect at one point $x$.
Then $x$ and $y = A^{-1}(x)$ satisfy
the desired property.
\end{proof}

Let $(\xi,\zeta,\mu,\nu)$ be as above.
We denote the corresponding submanifold of $Z_k(M)$ as
$Z_k(\xi,\zeta,\mu,\nu)$.
The intersection of two such submanifolds
can be completely determined as follows.

\begin{thm}
Two submanifolds $Z_k(\xi_i,\zeta_i,\mu_i,\nu_i)$, $i=1,2$
intersect if and only if
$$
      \frac{(\xi_1,\zeta_2)(\xi_2,\zeta_1)}{(\xi_1,\zeta_1)(\xi_2,\zeta_2)}
      = \frac{(\mu_1,\nu_2)(\mu_2,\nu_1)}{(\mu_1,\nu_1)(\mu_2,\nu_2)}.
 $$
When this condition is satisfied, the intersection is
one of the following:
\begin{align*}
	G_{k,n-1},& \quad \text{if $\xi_1=\xi_2$, $\mu_1=\mu_2$,}\\
	G_{k-1,n-1},& \quad \text{if $\zeta_1=\zeta_2$, $\nu_1=\nu_2$,}\\
	G_{k-1,n-2},& \quad \text{otherwise.}
\end{align*}
\end{thm}

\begin{proof}
We can put $\zeta_2 = x_0$ and $\nu_2 = y_0$.
By a similar argument in the proof of Theorem \ref{embd},
we can show that two submanifolds are disjoint if
$(\mu_1,y_0) = 0$.
Hence we can assume that $(\mu_1,y_0) \not= 0$.
Let $z$ and $w$ be elements of $G_k(x_0^\perp\oplus y_0^\perp)$
such that
$\cpsi_{a,b}(\xi_1,\zeta_1,\mu_1,\nu_1)(z) =
\cpsi_{a,b}(\xi_2, x_0, \mu_2, y_0)(w)$.
Then, since $\cphi_{0I,0J}(z) = 0$ for $\abs{I} = a, a-1$, we have $z = w$.
Let $\xi_1 = (1,\xi_1^1,\dots,\xi_1^l)$ and
$\mu_1 = (1,\mu_1^1,\dots,\mu_1^m)$. Then
the equations are written as
\begin{equation}\label{intcon}
\begin{split}
	\frac{\sum_{i\not\in I} {\bar\zeta}_1^i z^{iIJ'}}
	    { {\bar\zeta}_1^0 - \sum_i \xi_1^i {\bar\zeta}_1^i } +
	  \frac{\sum_{j\not\in J} {\bar\nu}_1^j z^{j'IJ'}}
	    { {\bar\nu}_1^0 + \sum_j \mu_1^j {\bar\nu}_1^j} &= 0,\\
	\sum_{i\in I} (\xi_1^i - \xi_2^i)z^{iIJ'} +
	  \sum_{j\in J} (\mu_1^j - \mu_2^j)z^{j'IJ'} &= 0.
\end{split}
\end{equation}
Let $L$ be the subspace of $x_0^\perp\oplus y_0^\perp$ corresponding
to $z$ and $\xi_i' = (0,\xi_i^1,\dots,\xi_i^l)$,
$\mu_i' = (0,\mu_i^1,\dots,\mu_i^m)$, $i=1,2$.
Then \eqref{intcon} means that
$$
	(\xi_1'-\xi_2', \mu_1'-\mu_2') \in L \subset
	(\frac{\zeta_1}{(\zeta_1, \xi_1)},
	   \frac{\nu_1}{(\nu_1, \mu_1)})^\perp.
 $$
Then the theorem follows immediately.
\end{proof}

\begin{cor}\label{cmlx}
 Let $M=\Hl\times \Proj^m$. Then we have
$$
	M_\C = \Hl\times {\bar\Hypb}^l \times
	  \left(\Proj^m\times \bar{\Proj}^m \setminus F(1,m;m+1)\right).
 $$
\end{cor}

\section{Some examples}
\begin{enumerate}
\item $M=\Proj(W)=\Proj^n$:
We can write $Z_k(\Proj^n) = F(k,k+1;n+1)$.
The projection map is written as
\begin{equation*}
\begin{split}
	F(k,k+1;n+1) &\rightarrow \Proj^n\\
	(U,U') &\mapsto U^\perp \cap U'
\end{split}
\end{equation*}
The complexification $M_\C$ is $\Proj^n \times \bar{\Proj}^n
\setminus F(1,n;n+1)$.
For $(\mu,\nu)\in M_\C$, the corresponding submanifold is written as
\begin{align*}\label{cmplx}
	\{(U, U') \mid U \subset \nu, \mu \subset U'\}.
\end{align*}
This is a straightforward generalization of the
well-known description of the twistor space of $\Proj^2$.

\item $M=\Hypb(V) = \Hypb^n$:
For a subspace $U$ of $V$, we write
$$
	U = \begin{cases}
	  U_{\pm},&\quad \text{if $U$ has a vector with positive
		norm},\\
	  U_{-},&\quad \text{if $G|_U$ is negative definite}.
	\end{cases}
 $$
Then we have
$$
	Z_k(\Hypb^n) =
	  \{ (U_-, U'_{\pm}) \in F(k,k+1;V) \}
 $$
The projection map is written as:
\begin{align*}
	Z_k(\Hypb^n) &\rightarrow \Hypb^n\\
	(U, U') &\mapsto U^\perp \cap U'
\end{align*}
$M_\C$ is $\Hypb^n \times \bar{\Hypb}^n$.
For $(\xi,\zeta)\in M_\C$, the corresponding submanifold
is written as:
$$
	\{(U,U')\mid U\subset \zeta, \xi \subset U' \}
 $$
\item $Z_1(\Hl \times \Proj^m)$ with $l,m\ge 1$:
The restriction of $\phi_{1,1}$ to $U_{1,0}\cap U_{0,1}$
is an isomorphism to
the projectivization of
the space of rank 2 elements of $V\otimes W$.
The boundary of the image is identified with
$\Proj(V)\times\Proj(W)$ by the Segre embedding.
The space $\Proj(V)$ has a decomposition:
$$
	\Proj(V) = \Hypb_+ \cup S^{2l-1} \cup \Hypb_-,
 $$
by the sign of the norm.
Then we obtain $Z_1(\Hl \times \Proj^m)$ by
changing the boundary as follows:
\begin{enumerate}
\item remove $S^{2l-1}\times \Proj^m$,
\item replace
  $\Hypb_+  \times \Proj^m$ with
    $\Proj(T^{(1,0)}\Hypb_+)\times \Proj^m$,
\item replace
$\Hypb_- \times \Proj^m$ with
$\Hypb_-\times \Proj(T^{(1,0)}\Proj^m)$.
\end{enumerate}
\end{enumerate}


\end{document}